\documentclass[a4paper,12pt]{article}
\usepackage[margin=3cm]{geometry}
\usepackage{amsmath}
\usepackage{amsthm}
\usepackage{amssymb}
\usepackage{color}
\usepackage{graphicx}

\def\Z{\mathbb Z}

\def\N{\mathbb N}
\def\A{\mathcal A}
\def\C{\mathcal C}
\def\L{\mathcal L}
\def\S{\mathcal S}
\def\pf{\begin{proof}}
\def\pfk{\end{proof}}
\newcommand{\Cd}{\vartriangle\!\C}

\newcommand{\lcp}[2]{\mathrm{lcp}(#1,#2)}

\newcommand{\vp}{\varphi}
\newcommand{\vpb}{\varphi_\beta}
\newcommand{\ubeta}{\mathbf{u}_\beta}
\newcommand{\Le}[1]{\mathrm{Lext}(#1)}
\newcommand{\Rex}[1]{\mathrm{Rext}(#1)}
\newcommand{\e}{\epsilon}
\newcommand{\ub}{\mathbf{u}}
\newcommand{\wb}{\mathbf{w}}

\newtheorem{prop}{Proposition}
\newtheorem{lem}[prop]{Lemma}
\newtheorem{rmrk}[prop]{Remark}
\newtheorem{thm}[prop]{Theorem}
\newtheorem{de}[prop]{Definition}
\newtheorem{coro}[prop]{Corollary}
\newtheorem{ass}[prop]{Assumption}

\begin{document}
\begin{center}
    \textbf{Factor complexity of infinite words associated with non-simple Parry numbers}
\end{center}

\begin{center}
    $\text{Karel Klouda}^{1,2}$, $\text{Edita
    Pelantov\'{a}}^{1,3}$ \\[2mm]
    \texttt{karel@kloudak.eu}
\end{center}

\small{$\ ^1$ FNSPE, Czech Technical University in Prague}

\small{$\ ^2$ LIAFA, Universit\'{e} Denis-Diderot (Paris VII)}

\small{$\ ^3$ Doppler Institute for mathematical physics and applied mathematics, Prague}\\[5mm]

The factor complexity of the infinite word $\ubeta$ canonically
associated to a non-simple Parry number $\beta$ is studied. Our
approach is based on the notion of special factors introduced by
Berstel and Cassaigne. At first, we give a handy method for
determining infinite left special branches; this method is
applicable to a broad class of infinite words which are fixed
points of a primitive substitution. In the second part of the
article, we focus on infinite words $\ubeta$ only. To complete the
description of its special factors, we define and study
$(a,b)$-maximal left special factors. This enables us to
characterize non-simple Parry numbers $\beta$ for which the word
$\ubeta$ has affine complexity.

\section{Introduction}

The aim of this work is to compute the \emph{factor complexity
function} $\C(n)$ of the infinite word $\ubeta$ associated with
$\beta$-expansions~\cite{Renyi1957}, where $\beta$ is a
\emph{non-simple Parry number}. The definition of Parry numbers is
connected with the R\'{e}nyi expansion of unity $d_\beta(1)$.
Parry numbers are those $\beta$ for which $d_\beta(1)$ is
eventually periodic. Positional numerical systems with a Parry
number as a base have a nice behavior. For example, if we consider
$\beta$-integers, i.e., real numbers with vanishing
$\beta$-fractional part in their $\beta$-expansion, then the
distances between two consecutive $\beta$-integers take only
finitely many values. In fact,  this property can be used as an
equivalent definition of Parry numbers. In this sense, positional
numeration systems based on Parry numbers are a natural
generalization of the classical decimal or binary systems. Let us
mention that even the innocent looking rational base $\beta =
\tfrac{3}{2}$ brings into numeration systems phenomena never
observed before~\cite{Akiyama2008}.

The most prominent Parry number is the golden mean $\tau =
\frac{1+\sqrt{5}}{2}$ with $d_\tau(1) =11$. The infinite word
associated to $\tau$ is the famous Fibonacci chain, i.e., the word
generated by the substitution $0\mapsto 01$ and $1\mapsto 0$. The
Fibonacci chain codes the distances between $\tau$-integers. Fabre
in~\cite{Fabre1995} showed that for any Parry number there exists
a canonical substitution over a finite alphabet such that its
unique fixed point $u_\beta$ represents the distribution of
$\beta$-integers on the real line.

$\beta$-integers attracted attention of physicists after the
discovery of quasicrystals in~1982~\cite{Shechtman1984}.
$\tau$-integers were shown to be a suitable tool for describing
atomic positions in solid materials with long range order and
non-crystalo\-graphical five-fold symmetry~\cite{Moody1993},
\cite{Barache1998}. The knowledge of the factor complexity of the
Fibonacci chain is the first step towards the description of
variability of local configurations in
quasicrystals~\cite{Masakova2003}.

Parry numbers are split into two groups: a Parry number $\beta$ is
called \emph{simple} if the R\'enyi expansion of unity
$d_\beta(1)$ has only a finite number of nonzero elements,
otherwise $\beta$ is \emph{non-simple}. The questions concerning
the factor complexity of words $\ubeta$ associated with simple
Parry numbers were discussed in \cite{Frougny2004} and
\cite{Bernat2007}. Of course, since among  $\ubeta$  one can find
some Sturmian sequences and Arnoux-Rauzy words, the complexity of
$\ubeta$ for some specific values of $\beta$ were known earlier.

The first non-simple Parry number $\beta$ for which the factor
complexity  of $\ubeta$ was precisely determined is such that
$d_\beta(1) = 2(01)^\omega$, i.e., $\beta$ is a root of the
polynomial $x^3-2x^2 - x +1$. This non-simple Parry number appears
naturally when describing the model of quasicrystals with
seven-fold symmetry~\cite{Frougny2003}. The first attempt to study
factor complexity of $\ubeta$ for broader class of non-simple
Parry number can be found in \cite{Frougny2007}.

Since any infinite word $\ubeta$ is the fixed point of a primitive
substitution, the factor complexity of $\ubeta$ can be estimate
from above by a linear function, see \cite{Queffelec1987}.
Moreover, we know that the first difference of complexity is
bounded by a constant~\cite{Cassaigne1996}~\cite{Mosse1996}.
Nevertheless, in general, it is hard to find an explicit formula
for the complexity function of an infinite word $\ub$ and it seems
it holds also for the case of $\ubeta$. However, we are able to
find all \emph{left special factors} that, in a certain sense,
completely determine the factor complexity. The notion of (right)
special factor was introduced by Berstel~\cite{Berstel1980} in
1980 and considerably enhanced by Cassaigne in his
paper~\cite{Cassaigne1997} in 1997. We introduce another slight
enhancement, a tool that will help us to identify all
\emph{infinite left special branches} of fixed point of
substitutions satisfying some natural assumption. Further, the
knowledge of the structure of left special factors will allow us
to identify all non-simple Parry numbers $\beta$ for which the
complexity of $\ubeta$ is affine: The complexity of $\ubeta$ is
affine if and only if $t_1(0\cdots0(t_1 - 1))^\omega$
(Theorem~\ref{thm:the_main_one}).

\section{Parry numbers and associated infinite words}

For each $x \in [0,1)$ and for each $\beta > 1$, using a greedy
algorithm, one can obtain the unique $\beta$-\emph{expansion}
$(x_i)_{i \geq 1}$,$x_i \in \N$, of the number $x$ such that
$$
    x = \sum_{i \geq 1} x_i \beta^{-i} \quad \text{and}\quad \sum_{i \geq k} x_i \beta^{-i} <
    \beta^{-k+1}.
$$
By shifting, each non-negative number has a $\beta$-expansion. For
$x \in [0,1)$, the $\beta$-expansion can be computed also by using
the piecewise linear map $T_\beta: [0,1) \rightarrow [0,1)$
defined as
$$
    T_\beta(x) = \{\beta x\},
$$
where $\{\beta x\}$ is the fractional part of the real number
$\beta x$. The sequence $d_\beta(x) = x_1 x_2 x_3 \cdots$ is
obtained by iterating $T_\beta$ with
$$
    x_i = \lfloor \beta T_\beta^{i-1}(x)\rfloor.
$$
The difference between $\beta$-expansion and $d_\beta(x)$ arises
for $x=1$ since the \emph{R\'enyi expansion of unity} $d_\beta(1)$
is not a $\beta$-expansion. Parry~\cite{Parry1960} showed that
$d_\beta(1)$ plays a very important role in the theory of
$\beta$-numeration. Among other things, it allows us to define
Parry numbers.
\begin{de}
    A real number $\beta > 1$ is said to be a Parry number if
    $d_\beta(1)$ is eventually periodic. In particular,
    \begin{itemize}
        \item[a)] if $d_\beta(1) = t_1 \cdots t_m$ is finite, i.e., it ends in infinitely many zeros, then
        $\beta$ is a \emph{simple Parry number},
        \item[b)] if it is not finite, i.e., $d_\beta(1) = t_1 \cdots t_m(t_{m+1}\cdots
        t_{m+p})^\omega$, then $\beta$ is called a \emph{non-simple Parry
        number}.
    \end{itemize}
\end{de}
Note, that the parameters $m,p > 0$ are taken the least possible.
It implies that $t_m \neq t_{m+p}$ which will be a very important
fact. Another crucial property of $d_\beta(1)$ is the following
\emph{Parry condition}~\cite{Parry1960} valid for all $\beta > 1$
\begin{equation}\label{eq:lex_order}
t_jt_{j+1}t_{j+2}\cdots \quad \prec \quad t_1 t_{2} t_{3}\cdots
\qquad\hbox{ for every}\ j> 1\,,
\end{equation}
where $\prec$ is the (strict) lexicographical ordering.

As the infinite word $\ubeta$ is tightly connected with a
geometrical interpretation of $\beta$-integers, we first introduce
$\beta$-integers along with some of their properties.
\begin{de}
    The real number $x$ is a \emph{$\beta$-integer} if the
    $\beta$-expansion of $|x|$ is of the form $\sum_{i=0}^k a_i \beta^i$. The set of all $\beta$-integers is denoted by
    $\Z_\beta$.
\end{de}
The definition of $\beta$-integers coincides with the definition
of classical integers in the case of $\beta$ in $\Z$. But there
are several new phenomena linked with the notion of
$\beta$-integers when $\beta$ is not an integer. For our purposes,
the most interesting difference between classical integers and
$\beta$-integers is the difference in their distribution on the
real line. While the classical integers are distributed
equidistantly, i.e., gaps between two consequent integers are
always of the same length 1, the lengths of gaps between
$\beta$-integers can take their values even in an infinite set.
More precisely, Thurston~\cite{Thurston1989} proved the following
theorem.
\begin{thm}
    Let $\beta > 1$ be a real number and $d_\beta(1) = (t_i)_{i \geq 1}$. Then the length of gaps between
    neighbors in $\Z_\beta$ takes values in the set $\{\triangle_0, \triangle_1,
    \ldots\}$, where
    $$
        \triangle_i = \sum_{k \geq 1} \frac{t_{k + i}}{\beta^k},
        \quad \text{for $i \in \N$.}
    $$
\end{thm}

\begin{coro}
    The set of lengths of gaps between two consecutive $\beta$-integers is
    finite if and only if $\beta$ is a Parry number.
    Moreover, if $\beta$ is a simple Parry number, i.e., $d_\beta(1) = t_1 \cdots t_m$, the set reads $\{\triangle_0, \triangle_1,
    \ldots \triangle_{m-1}\}$, if $\beta$ is a non-simple Parry
    number, i.e., $d_\beta(1) = t_1 \cdots t_m(t_{m+1}\cdots
        t_{m+p})^\omega$, we obtain $\{\triangle_0, \triangle_1,\ldots
    \triangle_{m + p-1}\}$.
\end{coro}
Now, let us suppose that we have drawn $\beta$-integers on the
real line and assume that $\beta$ is a Parry number. If we read
the length of gaps from zero to the right, we obtain an infinite
sequence, say $\{\triangle_{i_k}\}_{k \geq 0}$. Further, if we
read only indices, we obtain an infinite word over the alphabet
$\{0,\ldots,m-1\}$ in the case of simple Parry numbers, and over
the alphabet $\{0,\ldots,m+p-1\}$ in the non-simple case. The
obtained infinite word is just the word $\ubeta$ we are interested
in. However, there exists another way to define it.
Fabre~\cite{Fabre1995} proved that $\ubeta$ can be defined as the
unique fixed point of a substitution $\vp_\beta$ canonically
associated with a Parry number $\beta$ and defined as follows.
\begin{de} \label{def:vp_beta_simple}
For a simple Parry number $\beta$ the canonical substitution
$\varphi_\beta$  over the alphabet ${\cal A}=\{0,1,\dots,m-1\}$ is
defined by
            $$
            \begin{array}{ccl}
            \varphi_\beta(0)&=&0^{t_1}1\\
            \varphi_\beta(1)&=&0^{t_2}2\\
            &\vdots&\\
            \varphi_\beta(m\!-\!2)&=&0^{t_{m-1}}(m\!-\!1)\\
            \varphi_\beta(m\!-\!1)&=&0^{t_m}
            \end{array}
            $$
\end{de}
\begin{de} \label{def:vp_beta_nonsimple}
For a non-simple Parry number $\beta$ the canonical substitution
$\varphi_\beta$ over the alphabet ${\cal A}=\{0,1,\dots,m+p-1\}$
is defined by
            $$
            \begin{array}{ccl}
            \varphi_\beta(0)&=&0^{t_1}1\\
            \varphi_\beta(1)&=&0^{t_2}2\\
            &\vdots&\\
            \varphi_\beta(m\!-\!1)&=&0^{t_m}m\\
            \varphi_\beta(m)&=&0^{t_{m+1}}(m\!+\!1)\\
            &\vdots&\\
            \varphi_\beta(m\!+\!p\!-\!2)&=&0^{t_{m+p-1}}(m\!+\!p\!-\!1)\\
            \varphi_\beta(m\!+\!p\!-\!1)&=&0^{t_{m+p}}m.
            \end{array}
            $$
\end{de}
We see that the definition of $\vp_\beta$ is given by $d_\beta(1)$
and that the only difference between simple and non-simple cases
lies in the images of the last letters $m-1$ and $m+p-1$. While in
the simple case the last letters of images $\vp_\beta(k), k =
0,1,\ldots, m-1$, are all distinct and so the images form a
suffix-free code, in the non-simple case either $\vp_\beta(m) =
0^{t_m}m$ is a prefix of $\vp_\beta(m+p-1) = 0^{t_{m+p}}m$ or vice
versa. As we will see later on, this property is crucial from the
point of view of computing the complexity of the infinite
word~$\ubeta$.
\begin{de}
    Let $\beta > 1$ be a Parry number. The unique fixed point of
    the canonical substitution $\vp_\beta$ is denoted by
    $$
        \ubeta = \lim_{n \rightarrow \infty} \vp_\beta^n(0) =
        \vp^{\infty}_\beta(0).
    $$
\end{de}
The uniqueness of $\ubeta$ follows from the definitions of
$\vp_\beta$, the letter $0$ is the only admissible starting letter
of a fixed point.

\section{Special factors and factor complexity}

In this section, we will recall the notion of special factors of
an arbitrary infinite word and we will explain how the structure
of special factors of an infinite word determines its factor
complexity. To be able to do it, we need some usual basic
notation, see~\cite{Cassaigne1997} for more.
\begin{de}
    Let $\A = \{0,1,\ldots, q - 1\}, q \geq 1$ be a finite
    alphabet. An \emph{infinite word} over the alphabet $\A$ is a sequence
    $\ub = (u_i)_{i \geq 1}$ where $u_i \in \A$ for all $i \geq 1$. If
    $v = u_j u_{j+1}\cdots u_{j+n-1}$, $j,n \geq 1$, then $v$ is
    said to be a \emph{factor} of $\ub$ of length $n$ and the index $j$ is
    an \emph{occurrence} of $v$, the empty word $\e$ is the factor of length 0.

    By $\L_n(\ub)$ we denote the set of all factors of $\ub$ of length
    $n \in \N$, the \emph{language} of $\ub$ is then the set
    $\L(\ub)=\bigcup_{n\in\N}\L_n(\ub)$.
\end{de}

\begin{de}
    Let $\ub$ be an infinite word over an alphabet $\A$. The
    function $\C(n) = \#\L_n(\ub)$ is the \emph{factor complexity
    function} of $\ub$. We further define the first difference of the complexity by $\Cd(n) = \C(n+1) -
    \C(n)$.
\end{de}
In what follows, we shall restrict ourself to those infinite words
which are fixed point of some substitution (morphism) $\vp$
defined over a finite alphabet $\A$. We shall further assume that
$\vp$ is injective and \emph{primitive}.
\begin{de}
    A substitution $\vp$ is primitive if there exists $k \in \N$ such
    that for all $a,b \in \A$ the word $\vp^k(a)$ contains $b$.
\end{de}
Equivalently, $\vp$ is primitive if the incidence matrix $M_\vp$
is primitive.

There are several well-known properties of the complexity
function~$\C$.
\begin{prop} \label{prop:properties_of_C} \
    \begin{itemize}
        \item[(i)] For each infinite word $\ub$, $0 \leq \C(n) \leq
                    (\#\A)^n$,
        \item[(ii)] if $\ub$ is eventually periodic then $\C(n)$ is
                    eventually constant,
        \item[(iii)] $\ub$ is aperiodic in and only if $\C(n)$ is
                    unbounded and $\C(n)$ is unbounded if and only if $\Cd(n)  \geq 1$,
                    for all $n \in \N$,
        \item[(iv)] if $\ub$ is a fixed point of a primitive
                    substitution then $\C(n)$ is a sublinear
                    function, i.e., $\C(n) \leq an + b$, for some
                    $a,b \in \N$,
        \item[(v)]  if $\ub$ is a fixed point of primitive
                    substitution then $\Cd(n)$ is bounded.
    \end{itemize}
\end{prop}
Items $(i)-(iii)$ are obvious, $(iv)$ is due
to~\cite{Queffelec1987}, $(v)$ was proved in~\cite{Mosse1996} and
in a more general context in~\cite{Cassaigne1996}.

It is also well known that any fixed point of a primitive
substitution is uniformly recurrent, i.e., each factor occurs
infinitely many times and the gaps between its two consecutive
occurrences are bounded in length. It implies that each factor is
extendable both to the right and to the left.
\begin{de}
    Let $v$ be a factor of $\ub$, the set of \emph{left
    extensions} of $v$ is defined as
    $$
        \Le{v} = \{a \in \A \mid av \in \L{(\ub)}\}.
    $$
    If $\#\Le{v} \geq 2$, then $v$ is said to be a \emph{left
    special (LS)
    factor} of $\ub$.

    In the analogous way we define the set of \emph{right extensions}
    $\Rex{\ub}$ and a \emph{right special (RS) factor}.  If $v$ is both left
    and right special, then it is called \emph{bispecial}.
\end{de}
The connection between (left) special factors and the complexity
follows from the following reasoning. Let us suppose that
$\L_{n}(\ub) = \{v_1, \ldots, v_k\}, k \geq 1$ and let $\Le{v_i} =
\{a^{(i)}_{1}, \ldots, a^{(i)}_{\ell _i}\}$, $\ell_i \geq 1, i =
1,\ldots,k$. Now, it is not difficult to realize that
$$
\L_{n+1}(\ub) = \{a^{(1)}_{1}v_1, \ldots, a^{(1)}_{\ell
_1}v_1,a^{(2)}_{1}v_2,\ldots, a^{(k-1)}_{\ell
_{k-1}}v_{k-1},a^{(k)}_{1}v_k,\ldots,a^{(k)}_{\ell _{k}}v_k\},
$$
i.e., by concatenating all factors of length $n$ and all their
left extensions we obtain all factors of length $n+1$. It implies
that
\begin{equation} \label{eq:connection}
    \#\L_{n+1}(\ub) - \#\L_{n}(\ub) = \Cd(n) = \sum_{\begin{subarray}{c}v \in \L_n(\ub) \\ \text{$v$ is LS}\end{subarray}}
                \!\!\!(\#\Le{v} - 1).
\end{equation}
Hence, if we know all LS factors along with the number of their
left extensions, we are able to evaluate the complexity $\C(n)$
using this formula.

\subsection{Classification of LS factors}

Let $a,b \in \Le{v}$ be left extensions of a factor $v$ of $\ub$,
it means both $av$ and $bv$ are factors of $\ub$. If there exists
a letter $c \in \Rex{av} \cap \Rex{bv}$, we say that $v$ can be
extended to the right such that it remains LS with left extensions
$a,b$, indeed $a,b \in \Le{vc}$.
\begin{de}
    Let $a,b \in \Le{v}$ be distinct left extensions of a LS factor $v$ of
    $\ub$. We say that $v$ is an \emph{$(a,b)$-maximal LS factor} if $\Rex{av} \cap \Rex{bv} =
    \emptyset$, in words, $v$ cannot be extended to the right
    such that it remains LS with left extensions $a,b$.
\end{de}
\begin{figure}[h]
    \begin{center}\resizebox{12cm}{!}{\includegraphics{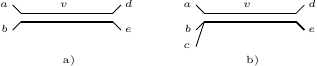}}\end{center}
    \caption{Two types of $(a,b)$-maximal LS factor $v$.} \label{mose_maximal}
\end{figure}
In general, there are two types of $(a,b)$-maximal LS factors both
depicted in Figure~\ref{mose_maximal}. In Case a{\small)}, $a$ and
$b$ are only left extensions of $v$ and so $v$ cannot be extended
to the right and remain LS. In Case b{\small)}, $v$ can be
prolonged by letter $e$ such that $ve$ is still a LS factor but it
looses its left extension $a$.

It can also happen that a factor $v$ with left extensions $a$ and
$b$ is extendable to the right infinitely many times. In this way
we obtain an infinite LS branch.
\begin{de}
    An infinite word $\mathbf{w}$ is an \emph{infinite LS
    branch} of $\ub$ if each prefix of $\mathbf{w}$ is~a~LS
    factor of $\ub$. We put
    $$
        \Le{\mathbf{w}} = \bigcap_{\text{$v$ prefix of $\mathbf{w}$}}\Le{v}.
    $$
\end{de}

\begin{prop} \
    \begin{itemize}
        \item[(i)]  If $\ub$ is eventually periodic, then there is
                    no infinite LS branch of $\ub$,
        \item[(ii)] if $\ub$ is aperiodic, then there exists at
        least one infinite LS branch of $\ub$,
        \item[(iii)]if $\ub$ is a fixed point of a primitive
                    substitution then the number of infinite LS
                    branches is bounded.
    \end{itemize}
\end{prop}
$(i)$ is obvious, $(iii)$ is a direct consequence
of~\eqref{eq:connection} and
Proposition~\ref{prop:properties_of_C} $(v)$. Item $(ii)$ is a
direct consequence of the famous K\"{o}nig's infinity
lemma~\cite{Konig1936} applied on sets $V_1, V_2, \ldots$, where
the set $V_k$ comprises all LS factors of length $k$ and where
$v_1 \in V_i$ is connected by an edge with $v_2 \in V_{i+1}$ if
$v_1$ is prefix of $v_2$.

Taking all together, our aim is to find all $(a,b)$-maximal LS
factors and also all infinite LS branches of $\ub$.

\begin{rmrk}
    The term ``special factor'' (for us it was RS factor) was introduced in 1980~\cite{Berstel1980} and it has been used for computing the factor complexity
    since then (eg. \cite{Berstel1989}, \cite{Luca1988}). The notations
    introduced above are based on Cassaigne's article~\cite{Cassaigne1997}. An $(a,b)$-maximal factor is a new term, actually it is a special case of a \emph{weak bispecial
    factor} proposed there. It is also shown in the article that bispecial factors determine the second difference of the
    complexity in a similar way as LS factors determine the first difference of the
    complexity.
\end{rmrk}

\begin{rmrk}
    Everything what has been (and will be) defined or showed for LS factors can be
    defined or showed similarly for RS factors.
\end{rmrk}

\subsection{How to find infinite LS branches}

\begin{figure}
    \begin{center}\resizebox{12cm}{!}{\includegraphics{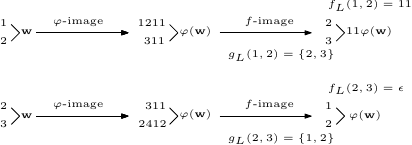}}\end{center}
    \caption{Images of LS factors.} \label{mose_f_image}
\end{figure}

Before introducing a new notion, let us consider the example
substitution
\begin{equation} \label{eq:ex_subst}
    \vp: 1 \mapsto 1211, 2 \mapsto 311, 3 \mapsto 2412, 4 \mapsto 435, 5 \mapsto
    534
\end{equation}
with $\ub = \vp^\infty(1)$. Further, let $w$ be a LS factor (or
infinite LS branch) of $\ub$ with left extensions $1$ and $2$. Is
$\vp(w)$ again LS factor? From Figure~\ref{mose_f_image} (first
line) we see that it is not since the letter~1 is its only left
extension. In order to obtain a LS factor, we have to prepend the
factor 11 which is the longest common suffix of $\vp(1) = 1211$
and $\vp(2) = 311$, then $11\vp(w)$ is a LS factor with left
extensions $2$ and $3$. In the case when $\Le{w} = \{2,3\}$
(second line in Figure~\ref{mose_f_image}), $\vp(w)$ is a LS
factor since the longest common suffix of $\vp(2) = 311$ and
$\vp(3) = 2412$ is the empty word $\e$.

\begin{de}
    Let $\vp$ be a substitution defined over an alphabet $\A$.
    For each couple of distinct letters $a,b \in \A$ we define
    $f_L(a,b)$ as the longest common suffix of words $\vp(a)$ and
    $\vp(b)$.
\end{de}

\begin{de}
   Let $v$ be a prefix of a word $w$, then $v^{-1}w$ is the word $w$ without the prefix
   $v$. Analogously, we define $wv^{-1}$, if $v$ is a suffix of $w$.
\end{de}

\begin{de} \label{def:g_L_ab_}
    Let $\vp$ be an injective substitution defined over an alphabet $\A$ having a fixed point $\ub$.
    For each unordered couple of distinct letters $a,b \in \A$ such that $\Rex{a} \cap \Rex{b} \neq
    \emptyset$ we define the set $g_L(a,b)$ as follows.
    \begin{itemize}
        \item[(i)]  If $f_L(a,b)$ is a proper suffix of both $\vp(a)$ and
                    $\vp(b)$, then $g_L(a,b)$ contains just the last letters of
                    factors $\vp(a)(f_L(a,b))^{-1}$ and $\vp(b)(f_L(a,b))^{-1}$.

        \item[(ii)] If $f_L(a,b) = \vp(a)$ (i.e., W.L.O.G. $|\vp(a)| < |\vp(b)|$),
                    then $g_L(a,b)$ contains the last letter of
                    the factor $\vp(b)(f_L(a,b))^{-1}$ and all the last letters of factors
                    $\vp(c)$, where $c \in \Le{a}$ such that $\Rex{ca} \cap \Rex{b} \neq
                    \emptyset$.
    \end{itemize}
\end{de}

\begin{ass} \label{ass:A}
    A substitution $\vp$ defined over $\A$ is injective and it has a fixed point
    $\ub$ such that for all $a,b \in \A$, for which $g_L(a,b)$ is
    defined, it holds that $\#g_L(a,b) = 2$.

    Moreover, if $f_L(a,b) = \vp(a)$ and $d$ is the last
    letter of the factor $\vp(b)(f_L(a,b))^{-1}$, then for all $c \in \Le{a}$
    such that $\Rex{ca} \cap \Rex{b} \neq \emptyset$ it holds that
    $d$ is not the last letter of $\vp(c)$.
\end{ass}
Assumption~\ref{ass:A} is valid for all suffix-free substitutions
since $g_L(a,b)$ from point $(i)$ of Definition~\ref{def:g_L_ab_}
contains always just two elements and the case when $f_L(a,b) =
\vp(a)$ never happens. If $f_L(a,b) = \vp(a)$, then
Assumption~\ref{ass:A} says that if $v$ is a LS factor with
$\Le{v} = \{a,b\}$, then the last letter $e$ of $\vp(c)$ is the
same for all $c \in \Le{av}$ and, moreover, $e\vp(a)$ is not a
suffix of $\vp(b)$ -- in other words, for each LS factor $v$ the
factor $f_L(a,b)\vp(v)$ is again LS. We will see that this
complicated assumption is satisfied for the (not suffix-free)
substitution $\vp_\beta$, where $\beta$ is a non-simple Parry
number.
\begin{de}
    Let $\vp$ be a substitution satisfying Assumption~\ref{ass:A}.
    Then for each LS factor (or infinite LS branch) $w$ having distinct left
    extensions $a$ and $b$ we define \emph{$f$-image} of $w$ as the
    factor $f_L(a,b)\vp(w)$.
\end{de}
With respect to the preceding discussion, Assumption~\ref{ass:A}
says that $f$-image is always LS factor and it has just two left
extensions, namely two elements of $g_L(a,b)$, corresponding to
two original left extensions $a$ and $b$.

\begin{figure}
    \begin{center}\resizebox{12cm}{!}{\includegraphics{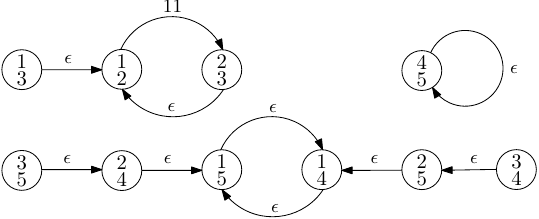}}\end{center}
    \caption{The graph $GL_\vp$ for the Substitution~\eqref{eq:ex_subst}.} \label{mose_GL_example}
\end{figure}

Assumption~\ref{ass:A} along with the notation introduced above
allow us to define the following graph.
\begin{de}
    Let $\vp$ be a substitution defined over an alphabet $\A$
    satisfying Assumption~\ref{ass:A}. We define a directed labelled graph
    $GL_\vp$as follows:

    \begin{itemize}
        \item[(i)] vertices of $GL_\vp$ are couples of distinct letters $a,b$ such that
            $\Rex{a} \cap \Rex{b} \neq \emptyset$,
        \item[(ii)] if $g_L(a,b) = \{c,d\}$, then there is an edge from a vertex $(a,b)$ to a vertex $(c,d)$ labelled by
        $f_L(a,b)$.
    \end{itemize}
\end{de}
In fact, the crucial result of Assumption~\ref{ass:A} is that
out-degree of each vertex is exactly one. The graph $GL_\vp$ for
our example substitution is drawn in Figure~\ref{mose_GL_example},
this substitution satisfies Assumption~\ref{ass:A} for it is
suffix-free.

Now, let us consider the case when $\wb$ is an infinite LS branch
with $a,b \in \Le{\wb}, a \neq b$. Obviously, $f$-image of $\wb$
is uniquely given. For most substitutions even a ``$f$-preimage''
of each infinite LS branch exists.
\begin{ass} \label{ass:B}
An infinite word $\ub$ is a fixed point of a substitution $\vp$
satisfying Assumption~\ref{ass:A}. For each infinite LS branch
$\wb$ of $\ub$ with $a,b \in \Le{\wb}, a \neq b$ there exists at
least one infinite LS branch $\overline{\wb}$ with left extensions
$c$ and $d$ such that $f$-image of $\overline{\wb}$ equals $\wb$
and $g_L(c,d) = \{a,b\}$.
\end{ass}
This assumption is very weak. Actually, we have not found any
primitive substitution not satisfying it. The reason for it is the
following. It is not difficult to prove (but it requires a lot of
new notation) this: Disruption of Assumption~\ref{ass:B} implies
that \emph{each} factor of $\ub$ can be decomposed into images of
letters $\vp(a), a \in \A$, in at least two different ways. Thus,
in order to prove that Assumption~\ref{ass:B} is satisfied, it
suffices to find any factor with unique decomposition to images of
letters. For the example substitution~\eqref{eq:ex_subst} the
factor 1211 is such a factor since it can arise only as $\vp(1)$,
an example of a factor having two different decomposition is 1112,
it is part of both $\vp(1)\vp(1)$ and $\vp(3)\vp(1)$. In the case
of $\vpb$, we will prove that Assumption~\ref{ass:B} is satisfied
differently.

\begin{thm}
    Let $\ub$ be a fixed point of a primitive injective
    substitution $\vp$ satisfying Assumption~\ref{ass:B} and let
    $\wb$ be an infinite LS branch with $a,b \in \Le{\wb}, a \neq
    b$. Then either $\wb$ is a \emph{periodic point} of $\vp$,
    i.e,
    \begin{equation}\label{eq:eq_for_inf_LS_type1}
        \wb = \vp^{\ell }(\wb) \quad \text{for some $\ell \geq 1$},
    \end{equation}
    and $(a,b)$ is a vertex of a cycle in $GL_\vp$ labelled by
    $\e$ only or $\wb = s\vp^{\ell }(s)\vp^{2\ell}(s)\cdots$ is the unique solution of the equation
    \begin{equation}\label{eq:eq_for_inf_LS_type2}
        \wb = s\vp^{\ell }(\wb),
    \end{equation}
    where $(a,b)$ is a vertex of a cycle in $GL_\vp$ containing at
    least one edge with non-empty label, $\ell$ is the length of this
    cycle and
    \begin{equation}\label{eq:def_of_prefix_s}
        s = f_L(g_L^{\ell -1}(a,b)) \cdots
        \vp^{\ell -2}(f_L(g_L(a,b))\vp^{\ell -1}(f_L(a,b)).
    \end{equation}
\end{thm}

\pf

\begin{figure}[h]
    \begin{center}\resizebox{7cm}{!}{\includegraphics{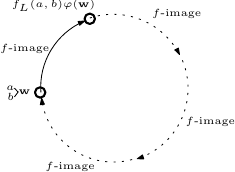}}\end{center}
    \caption{Circular structure of infinite LS branches.} \label{mose_structure_infinite}
\end{figure}

Due to Assumption~\ref{ass:B}, both the $f$-image and the
$f$-preimage of $\wb$ exist. The $f$-image is unique due to
Assumption~\ref{ass:A} and the uniqueness of $f$-preimage follows
from the fact that the number of infinite LS branches is finite.
Thus, $f$-image is one-to-one mapping on the finite set of all
ordered couples
$$
    \{((c,d),\overline{\wb})\},
$$
where $\overline{\wb}$ is an infinite LS branch of $\ub$ and
$(c,d)$ is an unordered couple of letters such that $c,d \in
\Le{\overline{\wb})}, c \neq d$. The $f$-image can be viewed as a
permutation on this finite set and so it decomposes the set to
independent cycles as depicted in
Figure~\ref{mose_structure_infinite}, $\ell$ is then the length of
the cycle containing $((a,b),\wb)$.

As explained earlier, applying $f$-image on a LS factor having
left extensions $a,b$ corresponds to the movement along the edge
in $GL_\vp$ which leads from $(a,b)$. If the labels of the edges
of the cycle are all $\e$, then the $f$-image coincides with $\vp$
and so we obtain the periodic
point~\eqref{eq:eq_for_inf_LS_type1}. If at least one edge is
labelled by non-empty word, then $\ell$ is also length of the
cycle in $GL_\vp$ containing the vertex $(a,b)$ and
Equation~\eqref{eq:eq_for_inf_LS_type2} then corresponds to
$\ell$-times applying of $f$-image on $\wb$.

\pfk

Our example substitution $\vp$ (see~\eqref{eq:ex_subst}) has five
periodic points
$$
\vp^\infty(1),\vp^\infty(4),\vp^\infty(5),(\vp^2)^\infty(2),(\vp^2)^\infty(3).
$$
It is an easy exercise to show that
\begin{multline*}
    \Le{1} = \{1,2,3,4,5\}, \Le{2} = \{1,4,5\}, \Le{3} = \{1,4,5\},\\
    \Le{4} = \{1,2,3\}, \Le{5} = \{1,2,3\}.
\end{multline*}
Looking at the graph $GL_\vp$ depicted in
Figure~\ref{mose_GL_example} we see that
$\vp^\infty(4),\vp^\infty(5)$ are not infinite LS branches as none
of the vertices $(1,2), (2,3)$ and $(1,3)$ is a vertex of a cycle
labelled by $\e$ only. Hence, only
$\vp^\infty(1),(\vp^2)^\infty(2),(\vp^2)^\infty(3)$ are infinite
LS branches with left extensions $1,4,5$.

As for infinite LS branches corresponding to
Equation~\eqref{eq:eq_for_inf_LS_type2}, in the case of our
example, there is only one cycle not labelled by the empty word
only between vertices $(1,2)$ and $(2,3)$. There are two (= the
length of the cycle) equations corresponding to this cycle
$$
    \wb = \vp(11)\vp^2(\wb) \quad \text{and} \quad \wb =
    11\vp^2(\wb).
$$
They give us two infinite LS branches
\begin{equation*}
    \begin{split}
        & \vp(11)\vp^3(11)\vp^5(11)\cdots,\\
        & 11\vp^2(11)\vp^4(11)\cdots,
\end{split}
\end{equation*}
the former having left extensions $1$ and $2$ and the latter $2$
and $3$.

\begin{rmrk}
    Assumption~\ref{ass:A} can be reformulated into a weaker form
    but to do so, it would require to introduce rather complicated
    notation. The important fact is that the canonical substitution
    $\vp_\beta$ satisfies Assumption~\ref{ass:A}.
\end{rmrk}

\section{Infinite LS branches of $\ubeta$}

\begin{figure}[h]
    \begin{center}\resizebox{7cm}{!}{\includegraphics{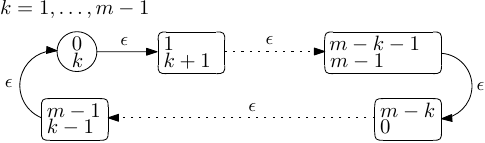}}\end{center}
    \caption{$GL_{\vp_\beta}$ for simple Parry $\beta$.} \label{mose_graph_simple}
\end{figure}

At first, let us recall known results for simple Parry numbers.
The substitution $\vp_\beta$ from
Definition~\ref{def:vp_beta_simple} is suffix-free and it implies
that it satisfies Assumption~\ref{ass:A}. One can easily prove
that even Assumption~\ref{ass:B} is satisfied. As mentioned
earlier, the last letters of images of letters are all distinct
and so $f_L(a,b) = \e$ for all couples $a,b \in \A$. The graph
$GL_{\vp_\beta}$ then looks as in Figure~\ref{mose_graph_simple}.
It contains $m-1$ cycles labelled by $\e$ only and hence the only
candidate for being an infinite LS branch is the unique fixed (and
periodic) point of $\vp_\beta$, namely $\ubeta$ with $\Le{\ubeta}
= \A$. The same result is proved in~\cite{Frougny2004} using
different techniques.

\begin{thm}[\cite{Bernat2007},\cite{Frougny2004}]
    Let $\beta > 1$ be a simple Parry number with $d_\beta(1) = t_1\cdots t_m$ and let $\ubeta$
    be the fixed point of the canonical substitution $\vp_\beta$~\eqref{def:vp_beta_simple}.
    Then
    \begin{itemize}
        \item[(i)] if $t_1=t_2=\cdots = t_{m-1} \quad \hbox{ or }\quad
        t_1>\max\{t_2,\dots,t_{m-1}\}$, the exact value of
        $\C(n)$ is known~\cite{Frougny2004},
        \item[(ii)] in particular, $(m-1)n+1\leq \C(n) \leq mn,\quad\hbox{ for all }\ n\geq
        1$,
        \item[(iii)] $\C(n)$ is affine $\Leftrightarrow$ the following two conditions are fulfilled \\
            {\rm 1)}\quad $t_m=1$ \\
            {\rm 2)}\quad for all $i=2,3,\ldots, m\!-\!1$ we have
            $$t_it_{i+1}\cdots t_{m-1}t_1\cdots t_{i-1} \quad \preceq \quad t_1t_2\cdots
            t_{m-1}.$$
        Then $\C(n)= (m-1)n + 1$.
    \end{itemize}
\end{thm}
In this paper, we will find the necessary and sufficient condition
for the complexity being affine in the case of non-simple Parry
numbers. We will see that it is more restrictive than the one from
point $(iii)$.

\subsection{Infinite LS branches in case of non-simple Parry numbers}

In this section, we will apply hitherto introduced theory on the
fixed point $\ubeta$ of the substitution $\vpb$, where $\beta$ is
a non-simple Parry number. To be able to do so, we need some more
notation and simple but useful technical lemmas.

\begin{de}
    For all $k, \ell \in \N$, we define the addition $\oplus : \N \times \N \rightarrow
    \A$ as follows.
    $$
        k \oplus \ell :=    \begin{cases}
                            k + \ell & \text{if $k + \ell < m + p$, } \\
                            m + (k + \ell - m \text{ mod } p) & \text{otherwise.}
                        \end{cases}
    $$
    Similarly, if used with parameters $t_i$, we define for all $k, \ell \in \N, k + \ell > 0$
    $$
        t_{k \oplus \ell} :=    \begin{cases}
                            t_{k + \ell} & \text{if $0 < k + \ell < m + p + 1$,} \\
                            t_{m + 1 + (k + \ell - m - 1 \text{ mod } p)} & \text{otherwise.}
                        \end{cases}
    $$
\end{de}
In fact, the addition $\oplus$ tracks the last letters of the
words $\vpb^n(0), n = 0,1,\ldots$. Therefore, we can rewrite the
definition of the substitution $\vp$ in a simpler form
$$
    \vpb(k) = 0^{t_{k + 1}}(k \oplus 1), \qquad \forall k \in \A.
$$
Further, employing the new notation and the definition of the
substitution $\vpb$, one can easily prove the following simple
observations.
\begin{lem} \label{lem:simple_observations}
For the substitution~$\vpb$ it holds
\begin{itemize}
    \item[(i)] for all $n \in \N$ and for all $k \in \A$
        \begin{equation*}
            \vpb^n(k) = (\vpb^{n-1}(0))^{t_{k \oplus 1}} (\vpb^{n-2}(0))^{t_{k
            \oplus 2}} \cdots (\vpb(0))^{t_{k \oplus (n - 1)}}0^{t_{k \oplus
            n}} (k \oplus n),
        \end{equation*}
    \item[(ii)] if $avb$ is a factor of $\ubeta$, $v \in \A^*$ and $a,b \neq 0$,
                then there exists unique factor $v'$ such that $\vpb(v') =
                vb$.
\end{itemize}
\end{lem}

Our aim is to obtain the graph $GL_{\vpb}$, thus, we need to know
left extensions of letters and also all $g_L(a,b)$.
\begin{de} \label{def:l_0_z(k)_and_co}
    Let us define for all $k \in \A$, $k \neq 0$, a function $z:\{1, \ldots, m+p-1\} \rightarrow
    \{0,1,\ldots,m+p-2\}$ by
    $$
        z(k) = \max\{j \in \N \mid 0^j \text{ is a suffix of } t_1
        t_2 \cdots t_k \}.
    $$
    For $k \in \{ m, \ldots, m+p-1 \}$ we also define a function $y:\{m, \ldots, m+p-1\} \rightarrow
    \{0,1,\ldots,p-1\}$ by
    $$
        y(k) = \begin{cases}
            \max\{j \in \N \mid 0^j \text{ is a suffix of } t_{m+1} t_{m+2} \cdots t_{m+p}t_{m+1} \cdots t_k\} & \text{if }k > m,\\
            \max\{j \in \N \mid 0^j \text{ is a suffix of } t_{m+1} t_{m+2} \cdots t_{m+p}\} & \text{if }k = m.
               \end{cases}
    $$
    Further, we define
    $$
        \ell_0 = \begin{cases}
                    0 & \text{if }t_1 > 1,\\
                    1 + \max\{j \in \N \mid 0^j \text{ is a prefix of } t_{2} t_{3} \cdots t_{m}\} & \text{otherwise}
               \end{cases}
    $$
    and finally we put $t = \min \{t_m,t_{m+p}\}$.
\end{de}
Note that $z(k)$ and $y(k)$ can return the same value for $k \geq
m$, a necessary condition for $z(k) \neq y(k)$ is that $t = 0$ and
$z(\ell) \neq y(\ell)$ for all $m \leq l < k$. Due to Parry
condition~\eqref{eq:lex_order} we must have $1\leq \ell_0 \leq
m-1$ as the case $d_\beta(1) = 10 \cdots 0(t_{m+1}\cdots
t_{m+p-1}1)^\omega$ is not admissible.
\begin{lem} \label{lem:lext_of_letters}
    For $\ubeta$ the fixed point of $\vpb$ it holds
    \begin{itemize}
    \item[(i)] $\Le{0} =  \{\ell _0, \ldots, m+p-1\}$,
    \item[(ii)] $\Le{k} = \{ z(k) \}$, for $k \in \{2,3, \ldots, m - 1 \}$,
    \item[(iii)] $\Le{k} = \{ z(k) , y(k) \}$, for $k \in \{m,m+1 \ldots, m + p - 1 \}$.
    \end{itemize}
\end{lem}

\pf

$(ii)$ Each letter $k > 0$ can appear in $\ubeta$ as the image of
$k-1$, namely $\vpb(k-1) = 0^{t_k}k$. If $t_k > 0$, then $0 \in
\Le{k}$, if $t_k = 0$ we consider $\vpb^2(k-2) =
\vpb(0^{t_{k-1}})k = (0^t_1 1)^{t_{k-1}}k$. Again, if $t_{k-1} >
0$, then $1 \in \Le{k}$, otherwise we continue in the same way.
Since $t_1 > 0$, this process is finite.

$(iii)$ The letter $m$ can appear in $\ubeta$ not only as an image
$\vpb(m-1)$ (i.e., case $(ii)$) but as well as $\vpb(m+p-1) =
0^{t_{m+p}}m$. If we realize this second possible origin of the
letters $m, m+1, \ldots, m+p-1$, then the proof is the same as for
$(ii)$.

$(i)$ If $t_1 > 1$, then $00$ is a factor of $\ub$. Hence, for all
$n \in \N$ the word $\vpb^n(0)0 = \cdots (0 \oplus n)0$ is a
factor as well and so $\Le{0} = \A$.

Let $t_1 = 1$, it implies $t_i \in \{0,1\}$ for $i =
1,\ldots,m+p$. It holds $\vpb^{\ell _0}(01) =
\vpb((\ell_0-1)\ell_0) = \ell_0 0 (\ell_0 + 1)$, hence, $\ell_0,
\ell_0+1, \ldots, m+p-1 \in \Le{0}$. But $d_ \beta(1)$ cannot
contain a sequence of consecutive 0's shorter than $\ell_0$ due to
Parry condition~\eqref{eq:lex_order} and so $\ell_0$ is the least
letter in $\Le{0}$.

\pfk

The previous lemma allows us to partially determine all LS factors
of $\ubeta$.
\begin{coro} \label{coro:beginnings_of_LS_factors}
    Let $v$ be a LS factor of $\ubeta$ containing at least one nonzero letter, then one of the following
    factors is a prefix of $v$.
    \begin{itemize}
      \item[(i)] $0^{t_1} 1$,
      \item[(ii)] $0^t m$,
      \item[(iii)] $0^{t_k}k$, if $k > m$ and $t = t_{m+1} = t_{m+2} = \cdots = t_{k-1} = 0$.
    \end{itemize}
\end{coro}
Note, that the factors from the last point are images of the
factor $0^t m$ in the case when $t=0$.
\begin{lem} \label{lem:g_L_f_L_for_vpb}
    For $\ubeta$ the fixed point of $\vpb$ it holds
    \begin{itemize}
      \item[(i)] if $(k,\ell)$ is unordered couple of distinct letters of $\A$ such that $\Rex{k} \cap \Rex{\ell } \neq \emptyset$,
                 and $(k,\ell) \neq (m-1,m+p-1)$, then $f_L(k,\ell) = \e$ and $g_L(k,\ell) = \{k \oplus 1, \ell \oplus 1\}$,
      \item[(ii)] $f_L(m-1,m+p-1) = 0^t m$ and $g_L(m-1,m+p-1) =  \{0,z\}$, where
                    \begin{equation}\label{eq:def_of_z}
                        z = \begin{cases}
                            1 + z(m-1)   & \text{if $t_{m} < t_{m+p}$}, \\
                            1 + z(m+p-1) & \text{if } t_{m+p} < t_{m}.
                        \end{cases}
                    \end{equation}

    \end{itemize}
\end{lem}

\pf

$(i)$ follows directly from the definitions of $g_L, f_L$ and
$\vpb$. $(ii)$ is a simple consequence of
Lemma~\ref{lem:lext_of_letters}. Remark that if $t_m > t_{m+p}
\geq 0$, then $z(m+p-1) = y(m+p-1)$.

\pfk

Now we have the knowledge necessary to complete the graph
$GL_{\vpb}$ but still we have to prove that the substitution
$\vpb$ satisfies Assumptions~\ref{ass:A} and~\ref{ass:B}.
\begin{lem}
    The substitution $\vpb$ from
    Definition~\ref{def:vp_beta_nonsimple} satisfies Assumptions~\ref{ass:A} and~\ref{ass:B}.
\end{lem}

\pf
    The fact that Assumptions~\ref{ass:A} is fulfilled follows from
    Lemmas~\ref{lem:lext_of_letters} and~\ref{lem:g_L_f_L_for_vpb}.

    To construct an $f$-preimage for an arbitrary infinite LS
    branch is easy due to Lemma~\ref{lem:simple_observations} part
    $(ii)$, Corollary~\ref{coro:beginnings_of_LS_factors} and
    Lemma~\ref{lem:g_L_f_L_for_vpb}.

\pfk

\begin{figure}[h]
    \begin{center}\resizebox{10cm}{!}{\includegraphics{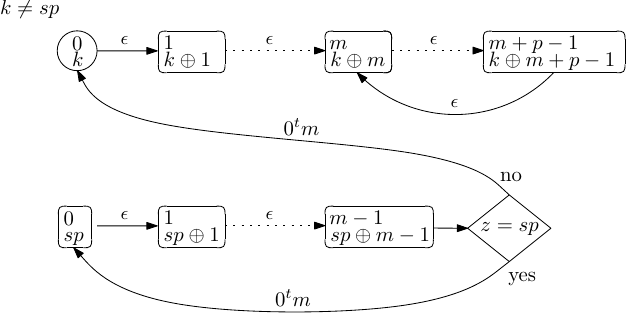}}\end{center}
    \caption{$GL_{\vp_\beta}$ for non-simple Parry $\beta$, $s$ is a positive integer.} \label{mose_graph_nonsimple}
\end{figure}

Now, we know all we need to be able to construct the graph
$GL_{\vpb}$. For the case when $t_1 > 1$, the graph is depicted in
Figure~\ref{mose_graph_nonsimple}. Since $\Le{0} = \A$, all
possible unordered couples of letters are vertices of the graph.
If $z$ is not a multiple of $p$ (i.e., the decision condition $z =
sp$ in Figure~\ref{mose_graph_nonsimple} returns \emph{no}), then
the graph contains only cycles with edges labelled by $\e$ only.
If $z = sp$ for certain positive integer $s$, then there is the
cycle on vertices $(0,z), (1,z \oplus 1), \ldots, (m-1, z \oplus m
- 1)$, where the edge from the vertex $(m-1, z \oplus m - 1)$ to
the vertex $g_L(m-1, z \oplus m - 1) = (0,z)$ is labelled by
$f_L(m-1, z \oplus m - 1) = 0^tm$.

If $t_1 = 1$ the graph $GL_{\vpb}$ is the same as in
Figure~\ref{mose_graph_nonsimple} but we have to remove vertices
$(k,\ell)$, where $k < \ell_0$ or $\ell < \ell_0$ and $(k,\ell)
\neq (0 \oplus n, z \oplus n)$ for any $n \in \N$. What is
important for our purpose is that the structure of cycles is the
same for arbitrary value of $t_1$.

Since the fact whether $z$ is or is not a multiple of $p$ is
crucial for the structure of cycles in $GL_{\vpb}$, we introduce
the following set.
\begin{de} \label{def:set_S}
    A non-simple Parry number $\beta > 1$ is an element of a set $\S$
    if and only if there exists a positive integer $s$ such that $z = sp$,
    where $z$ is the non-zero left extension of $0^t m$.
\end{de}
Employing Lemmas~\ref{lem:lext_of_letters}
and~\ref{lem:g_L_f_L_for_vpb}, one can easily prove the following.
\begin{lem}
    A non-simple Parry number $\beta > 1$ belongs to
    $\mathcal{S}$ if and only if one of the following conditions is satisfied
    \begin{equation*}
        \begin{split}
            \text{a{\small)} } & d_{\beta}(1) = t_1\cdots t_m(0\cdots0t_{m+p})^{\omega} \quad \text{and $t_m > t_{m+p}$}, \\
            \text{b{\small)} } & d_{\beta}(1) = t_1\cdots \underbrace{t_{m-qp}}_{\neq 0} \underbrace{0\cdots0}_{qp-1} t_m(t_m+1\cdots t_{m+p})^{\omega},\quad \text{$q \geq 1$, and $t_m <
            t_{m+p}$}.
        \end{split}
    \end{equation*}
\end{lem}
Putting it all together, we obtain a proof of the following
proposition which gives us the complete list of infinite LS
branches of $\ubeta$ for all non-simple Parry numbers.
\begin{prop} \label{prop:list_of_inf_branches}
    Let $\beta > 1$ be a non-simple Parry number and let $\ubeta$
    be the fixed point of the canonical substitution $\vp_\beta$.
    Then
    \begin{itemize}
        \item[(i)] if $p > 1$, then $\ubeta$ is an infinite LS
        branch with left extensions $\{m, m+1, \ldots,m+p-1\}$,
        \item[(ii)] if $\beta \notin \mathcal{S}$, then
        $\ubeta$ is the unique infinite LS branch,
        \item[(iii)] if $\beta \in \mathcal{S}$, then
        there are $m$ infinite LS branches
        \begin{equation*}\label{eq:branches}
            \begin{split}
                & 0^t m \vp^m(0^t m) \vp^{2m}(0^t m) \cdots \\
                & \quad \quad \quad  \quad \quad \vdots \\
                & \vp^{m-1}(0^t m) \vp^{2m-1}(0^t m) \vp^{3m-1}(0^t m)
                \cdots.
            \end{split}
        \end{equation*}
        There are no other infinite LS branches of $\ubeta$.
    \end{itemize}
\end{prop}

\section{Maximal LS factors}

As explained earlier, in order to determine the complexity of an
infinite word, we need to find all infinite LS branches as well as
all $(a,b)$-maximal LS factors. The structure of $(a,b)$-maximal
LS factors is not so simple as the one of infinite LS branches but
still it can be described using the notion of $f$-image. To define
an $f$-image for $(a,b)$-maximal LS factors, we need
Assumption~\ref{ass:A} to be satisfied also for $g_R$ -- we will
say that the \emph{right version of} Assumption~\ref{ass:A} is
satisfied.
\begin{lem}
    For the substitution $\vpb$ and for all distinct $a,b \in \A$ we
    have $f_R(a,b) = 0^{t_{a,b}}$, where
    \begin{equation}\label{eq:def_of_t_ab}
        t_{a,b} = \min\{t_a,t_b\}.
    \end{equation}
\end{lem}
Thus, the right version of Assumption~\ref{ass:A} is satisfied for
$\vpb$ is prefix-free.
\begin{de}
    A factor $v \in \A^+$ is an \emph{$(a-c,b-d)$-bispecial} factor of an
    infinite word $\ub$ defined over a finite alphabet $\A$ if
    both $avc$ and $bvd$ are factors of $\ub$.
\end{de}

\begin{de}
    Let a substitution $\vp$ defined over a finite alphabet $\A$ satisfy the left and right version of
    Assumption~\ref{ass:A} and let $v$ be an $(a-c,b-d)$-bispecial
    factor of a fixed point of $\vp$. Then
    $f_L(a,b)\vp(v)f_R(c,d)$ is said to be the \emph{$f$-image} of $v$.
\end{de}
Obviously, the $f$-image of $v$ is
$(\tilde{a}-\tilde{c},\tilde{b}-\tilde{d})$-bispecial, where
$g_L(a,b) = \{\tilde{a},\tilde{b}\}$ and $g_R(c,d) =
\{\tilde{c},\tilde{d}\}$.

Now, consider again the particular case of $\ubeta$. A LS factor
$v$ having $a,b \in \Le{v}$ is $(a,b)$-maximal if $\Rex{av} \cap
\Rex{bv} = \emptyset$ and so it is as well an
$(a-c,b-d)$-bispecial for all $c \in \Rex{av}$ and $d \in
\Rex{bv}$. Are $f$-images of $v$ again $(g_L(a,b))$-maximal? Not
all of them as states the following simple lemma.
\begin{lem} \label{lem:max_f_images_of_bispecials}
Let $v$ be a bispecial factor of $\ubeta$ having left extensions
$a$ and $b$. If its $f$-image
$$
    f_L(a,b)\vpb(v)f_R(c,d) = f_L(a,b)\vpb(v)0^{t_{c\oplus1,d\oplus1}},
$$
is $(g_L(a,b))$-maximal, then $c \in \Rex{av}, d \in \Rex{bv}$
satisfy
\begin{equation}\label{eq:max_f_image}
    \begin{split}
        t_{c \oplus 1} & \geq \max\{t_{e \oplus 1,f \oplus 1} \mid e \in \Rex{av}, f \in \Rex{bv}\} \\
        t_{d \oplus 1} & \geq \max\{t_{e \oplus 1,f \oplus 1} \mid e \in \Rex{bv}, f \in \Rex{bv}\}.
    \end{split}
\end{equation}
\end{lem}
%
%
\begin{de}
    An $f$-image of a bispecial factor $v$ having left extensions $a$ and $b$
    $$
        f_L(a,b)\vpb(v)f_R(c,d),
    $$
    where $c \in \Rex{av}, d \in \Rex{bv}$ satisfy~\eqref{eq:max_f_image}, is said to be the \emph{max-$f$-image} of $v$.
\end{de}
The following lemma is crucial for understanding the structure of
the max-$f$-images of $(a,b)$-maximal factors.
\begin{lem}\label{lem:main_lemma}
    If $\ell,k \in \A$, $\ell \neq k$ and $t_{\ell  \oplus 1} t_{\ell  \oplus 2}\cdots \succeq t_{k \oplus 1} t_{k \oplus
    2}\cdots$, then for all $n \in \N$ the longest common prefix
    of the factors $\vpb^{n}(k)$ and $\vpb^{n}(\ell)$, denoted by $\lcp{\vpb^{n}(k)}{\vpb^{n}(\ell)}$, equals
    \begin{equation*}
        \lcp{\vpb^{n}(k)}{\vpb^{n}(\ell)} = \vpb^{n}(k)(k \oplus
        n)^{-1},
    \end{equation*}
    i.e., $\vpb^{n}(k)$ without the last letter $k \oplus n$.

    Moreover, denote by $c$ the letter such that $(\lcp{\vpb^{n}(k)}{\vpb^{n}(\ell)})c$ is a prefix of $\vpb^n(\ell)$.
    Then, $t_{c \oplus 1} t_{c \oplus 2}\cdots \succeq t_{k \oplus (n+1)} t_{k \oplus
    (n+2)}\cdots$ for all $n \in \N$.
\end{lem}

\pf

The case $n = 0$ is trivial.

The rest of the proof is carried on by induction on $n$.
\begin{eqnarray}
  \vpb^{n+1}(k) &=& (\vpb^{n}(0))^{t_{k \oplus 1}}\vpb^{n}({k \oplus 1}), \nonumber \\
  \vpb^{n+1}(\ell) &=& (\vpb^{n}(0))^{t_{k \oplus 1}}(\vpb^{n}(0))^{t_{\ell  \oplus 1} -
  t_{k \oplus 1}}\vpb^{n}({\ell  \oplus 1}), \label{eq:lext_of_lcp}
\end{eqnarray}
if $t_{\ell  \oplus 1} = t_{k \oplus 1}$, we apply the assumption
of induction on $\lcp{\vpb^{n}({k \oplus 1})}{\vpb^{n}({\ell
\oplus 1})}$ and if $t_{\ell  \oplus 1} > t_{k \oplus 1}$, then on
$\lcp{\vpb^{n}({k \oplus 1})}{\vpb^{n}(0)}$ (see Parry
condition~\eqref{eq:lex_order}).

As for the second part of the statement, the letter $c$ is given
by~\eqref{eq:lext_of_lcp} and this along with the Parry condition
concludes the proof.

\pfk

\begin{lem} \label{lem:max_f_images_of_bispecials_explicit}
    Let $n \in \N$. The $n$-th max-$f$-image of a bispecial factor $v$ with left extensions $a$ and $b$,
    i.e., the factor we obtain if we apply $n$ times the mapping max-$f$-image on $v$,  equals
    $$
       \overline{v} = s\vpb^n(v)\lcp{\vpb^{n}(c)}{\vpb^{n}(d)},
    $$
    where $c \in \Rex{av}$, $d \in \Rex{bv}$, $s$ is given by (cf.~\eqref{eq:def_of_prefix_s})
    \begin{equation}\label{eq:def_of_prefix_s_with_n}
        s = f_L(g_L^{n-1}(a,b)) \cdots
        \vp^{n-2}(f_L(g_L(a,b))\vp^{n-1}(f_L(a,b)).
    \end{equation}
    and
    \begin{eqnarray*}
      t_{c \oplus 1} t_{c \oplus 2} \cdots & \succeq & t_{c' \oplus 1} t_{c' \oplus 2}\cdots, \\
      t_{d \oplus 1} t_{d \oplus 2} \cdots & \succeq & t_{d' \oplus 1} t_{d' \oplus 2}\cdots
    \end{eqnarray*}
    for all $c' \in \Rex{av}$ and $d' \in \Rex{bv}$.
\end{lem}

\pf

The case $n = 0$ is obvious, we carry on by induction on $n$. Let
us assume W.L.O.G. that
$$
    t_{c \oplus 1} t_{c \oplus 2} \cdots  \succeq t_{d \oplus 1} t_{d \oplus 2}\cdots
$$
and that $g_L^n(a,b) = \{\tilde{a},\tilde{b}\}$. Hence
$$
    \overline{v} = s\vpb^n(v)\vpb^{n}(d)(d \oplus n)^{-1}
$$
and
$$
\Rex{\tilde{b}\overline{v}} = \{d' \oplus n \mid t_{d' \oplus 1}
\cdots t_{d' \oplus n} = t_{d \oplus 1} \cdots t_{d \oplus n}\}.
$$
Further, let $c' \in \Rex{\tilde{a}\overline{v}}$, then due to
Lemma~\ref{lem:main_lemma}
$$
    t_{c' \oplus 1} t_{c' \oplus 2} \cdots  \succeq t_{d' \oplus (n + 1)} t_{d' \oplus (n +
    2)}\cdots
$$
for all $d' \oplus n \in \Rex{\tilde{b}\overline{v}}$. But $t_{d
\oplus (n + 1)} \geq t_{d' \oplus (n + 1)}$ for all $d' \oplus n
\in \Rex{\tilde{b}\overline{v}}$ and so the max-$f$-image of
$\overline{v}$ equals
$$
    f_L(g_L^n(a,b))\vpb(\overline{v})0^{t_{d \oplus (n + 1)}} =
    f_L(g_L^n(a,b))\vpb(s)\vpb^{n+1}(v)\lcp{\vpb^{n+1}(c)}{\vpb^{n+1}(d)}.
$$

\pfk

Each bispecial factor $v$ having left extensions $a$ and $b$ has
the unique max-$f$-image. Since the substitution $\vpb$ is
injective, the structure of max-$f$-images cannot be circular as
it is for $f$-images of infinite LS branches -- $v$ cannot be the
$k$-th max-$f$-image of its own for any $k$. However, the notion
of a max-$f$-image allows us to describe all $(a,b)$-maximal
factors of $\ubeta$ for all $a,b \in \A$. We will prove that each
$(a,b)$-maximal factor is the $k$-th max-$f$-image either of
$0^{t_1-1}$ if $t_1 > 1$ or of $0$ if $t_1=1$, for some $k \in
\N$. A sketch of the proof is as follows. Let $v$ be an
$(a,b)$-maximal factor containing at least two nonzero letters.
Employing Lemma~\ref{lem:simple_observations} part $(ii)$, one can
find a bispecial factor $\overline{v}$ such that its max-$f$-image
is $v$. Again, if $\overline{v}$ contains at least two nonzero
letters, we find a bispecial factor $\overline{\overline{v}}$ such
that its max-$f$-image is $\overline{v}$. In this way, we obtain a
bispecial factor containing at most one nonzero letter such that
its $k$-th max-$f$-image equals $v$. According to
Corollary~\ref{coro:beginnings_of_LS_factors}, the only candidates
for such bispecials are of the form $0^s$ or $0^tm0^q$, where $1
\leq s \leq t_1$ and $0 \leq q \leq t_1$. Note that $0^{t_1 + 1}$
cannot be a factor of $\ubeta$ and that is why we consider $s,q
\leq t_1$. In the case when $t = 0$, $0^{t_k}k0^q, k > m, t_{m+1}
= \cdots = t_{k-1} = 0$, could also be taken as the candidates but
we do not consider them as they are just prefixes of
$\vpb^{k-m}(m0^q)$. The following two lemmas tell us that $0^{t_1
- 1}$ (resp. $0$ if $t_1 = 1$) is the only candidate.
\begin{lem}
    Let $t_1 > 1$ and $k \in \N$. Then the $k$-th max-$f$-image of factors $0^{t_1}$, $0^s$ and $0^tm0^q$, where $1
    \leq s < t_1 - 1$ and $0 \leq q \leq t_1$, are not $(a,b)$-maximal
    for any distinct letters $a$ and $b$.
\end{lem}

\pf

First, consider $0^{t_1}$ with distinct left extensions $a$ and
$b$. It holds that $\Le{0^{t_1}} = \Le{0^{t_1}1}$ and
$\Rex{0^{t_1}} \subset \{k \in \A \setminus \{0\} \mid t_k =
t_1\}$. For each $k \in \Rex{0^{t_1}}$, we must have $t_{k \oplus
1} t_{k \oplus 2} \cdots \prec t_2t_3\cdots$ (see Parry
condition~\eqref{eq:lex_order}) and, due to
Lemma~\ref{lem:max_f_images_of_bispecials_explicit}, $k$-th
max-$f$-image of $0^{t_1}$ is a prefix of the $k$-th $f$-image of
the LS factor $\vpb^k(0^{t_1}1)$, both having the same left
extensions.

Similar arguments can be used in order to prove that $k$-th
max-$f$-image of $0^s$ is always a prefix of $k$-th $f$-image of
the LS factor $0^{t_1 - 1}$. Again, $\Le{0^s} = \Le{0^{t_1-1}}$
and the rest is implied directly by the Parry condition.

Finally, consider the LS factor $0^tm0^q$ having just two left
extensions $0$ and $z$ (see~\eqref{eq:def_of_z}). In accord with
Lemma~\ref{lem:max_f_images_of_bispecials_explicit}, the $m$-th
max-$f$-image of $0^{t_1 - 1}$ with left extensions $0$ and $p$
equals
\begin{equation}\label{eq:the_max_factor}
    0^t m \vpb^m(0^{t_1 - 1})\vpb^m(1)(m+1)^{-1} = 0^t m 0^{t_1} 1
    \cdots.
\end{equation}
Indeed, $\Rex{00^{t_1-1}} = \{k \in \A \setminus \{0\} | t_k = t_1
\}$ and $0 \in \Rex{p0^{t_1 -1}}$ and so the fact that $t_{k
\oplus 1} t_{k \oplus 2} \cdots \prec t_2t_3\cdots$ and the Parry
condition imply that the $m$-th max-$f$-image is $0^tm
\vpb^m(0^{t_1 - 1}) \lcp{\vpb^m(0)}{\vpb^m(1)}$. Thus, $0^tm0^q$,
as a prefix of~\eqref{eq:the_max_factor}, is not a
$(0,z)$-maximal.

\pfk

\begin{lem}
    Let $t_1 = 1$ and $k \in \N$. Then $t = 0$ and the $k$-th max-$f$-image of the factor $m0^q$,
    where $0 \leq q \leq 1$, is not $(a,b)$-maximal
    for any distinct letters $a$ and $b$.
\end{lem}

\pf

As in the proof of the previous lemma, we can prove that the
$(m-\ell_0)$-th max-$f$-image of $0$ with left extensions $\ell_0$
and $\ell_0 + p$ is the factor
\begin{equation}\label{eq:the_max_factor_for_t_0}
    m\vpb^m(1)(m+1)^{-1},
\end{equation}
where, according to Lemma~\ref{lem:simple_observations}
part~$(i)$,
$$
    \vpb^m(1) = (\vpb^{m-1}(0))^{t_{2}} (\vpb^{m-2}(0))^{t_{3}} \cdots
                (\vpb(0))^{t_{m}}0^{t_{m+1}}(m+1).
$$
In order that $m0$ may be $(0,z)$-maximal, it must be $\vpb^m(1) =
m+1$ and so $t_2 = \cdots = t_{m+1} = 0$. But is is not possible
due to the Parry condition since then $t_1 t_2\cdots  \cdots \prec
(t_{m+p}=1)t_m \cdots t_{m+p}t_{(m+p) \oplus 1}\cdots$.

\pfk

\begin{prop} \label{lem:generators_of_maximals}
    Let $v$ be an $(a,b)$-maximal factor of $\ubeta$. Then there
    exists $k \in \N$ such that $v$ is the $k$-th max-$f$-image of
    \begin{itemize}
        \item[(i)] $0^{t_1 - 1}$ if $t_1 > 1$,

        \item[(ii)] $0$ if $t_1 = 1$.
    \end{itemize}
\end{prop}

\pf

We will prove that if $v$ contains at least two nonzero letters,
then it is the $k$-th max-$f$-image of a bispecial factor of the
form $0^s$ or $0^tm0^q$, where $1 \leq s \leq t_1$ and $0 \leq q
\leq t_1$. The rest of the proof then follows from the previous
two lemmas.

Let us assume that $v$ contains at least two nonzero letters.
Then, due to Lemma~\ref{lem:simple_observations} part $(ii)$, $v =
f_L(a',b')\vpb(\overline{v})f_R(c',d')$, where $\overline{v}$ is a
$(a'-c',b'-d')$-bispecial factor such that $v$ is the
max-$f$-image of $\overline{v}$ and $g_L(a',b') = \{a,b\}$.
Analogously, if $\overline{v}$ contains at least two nonzero
letters, there exists an $(a''-c'',b''-d'')$-bispecial factor
$\overline{\overline{v}}$ which is an $f$-preimage of
$\overline{v}$. But it must be also a max-$f$-preimage, if it is
not, then $\overline{v}0^{q'}$ is also the $f$-image of
$\overline{\overline{v}}$ having the left extensions $a'$ and $b'$
for some $q' > 0$ and so $v$ cannot be $(a,b)$-maximal as it is a
proper prefix of the max-$f$-image of LS factor
$\overline{v}0^{q'}$ with the left extensions $a$ and $b$. Using
this argument iteratively, we will obtain a bispecial factor of
the form $0^s$ or $0^tm0^q$ such that $v$ is its $k$-th
max-$f$-image.

\pfk

In fact, the previous proposition along with
Lemma~\ref{lem:max_f_images_of_bispecials_explicit} provides us
with the complete list of $(a,b)$-maximal factors. However, in the
last section of this paper we will need to know some details to be
able to determine under which conditions the complexity of
$\ubeta$ is affine.
\begin{coro} \label{prop:maximals_0t_m}
If $d_\beta(1) \neq t_1(0\cdots0(t_1 - 1))^\omega$, then the
$k$-th max-$f$-image of the factor~\eqref{eq:the_max_factor} is
$(g_L^k(0,z))$-maximal for all $k \in \N$.

If $\beta \notin \S$, then the $k$-th max-$f$-image reads
$$
    \vpb^k(0^t m) \vpb^{m+1}(0^{t_1 - 1})\vpb^{m+1}(1)(m \oplus k)^{-1}.
$$
\end{coro}

\pf

The factor~\eqref{eq:the_max_factor} is always LS with just two
left extensions $0$ and $z$. Therefore it is $(0,z)$-maximal if it
is neither a prefix of any infinite LS branch or a proper prefix
of the $k$-th max-$f$-image of its own for certain $k > 0$.

In the case when $\beta \notin \S$, the longest common prefix of
the $k$-th max-$f$-image of the factor~\eqref{eq:the_max_factor}
and of the unique infinite LS branch $\ubeta$ equals
$$
    \vpb^k(0^t m)(m \oplus k)^{-1}.
$$
Hence, either it is non-empty and shorter than the longest common
prefix of the $(k+1)$-th max-$f$-image
of~\eqref{eq:the_max_factor} and of $\ubeta$ or it is empty, $k <
p$ and $t=t_{m+1}=\cdots=t_{m+k}=0$ (or only $t = 0$ for $k=0$).
In the latter case, the $k$-th max-$f$-image
of~\eqref{eq:the_max_factor} begins in letter $m+k$ which is
different from the first letters of $\ubeta$ and of all other
max-$f$-images of~\eqref{eq:the_max_factor}. Putting all together,
the $k$-th max-$f$-image of~\eqref{eq:the_max_factor} is neither a
prefix of $\ubeta$ or of the $\ell$-th max-$f$-image
of~\eqref{eq:the_max_factor} for any $\ell \neq k$.

If $\beta \in \S$, then $\ubeta$ is not the only one infinite LS
branch, there are $m$ other branches
\begin{equation}\label{eq:ub_1}
\ub_1 = 0^t m \vpb^m(0^t m) \vpb^{2m}(0^t m) \cdots
\end{equation}
and $\ub_\ell = \vpb^{\ell -1}(\ub_1), l = 2,\ldots,m$. To finish
the proof, we have to foreclose the possibility that the
factor~\eqref{eq:the_max_factor} is prefix of $\ub_1$. Looking
at~\eqref{eq:ub_1} and~\eqref{eq:the_max_factor}, we see that it
happens only if $t=t_1-1$ and $m=1$, in other words, if
$d_\beta(1) = t_1(0\cdots0(t_1 - 1))^\omega$. The proof of that
the factor~\eqref{eq:the_max_factor} is not a prefix of any
max-$f$-image of its own is analogous to the one above.

\pfk

\begin{coro} \label{prop:maximals_t_1_gr_1}
If $t_1 > 1$, then the $k$-th max-$f$-image of $0^{t_1-1}$ with
left extensions $0$ and $a$ is a $(g_L^k(0,a))$-maximal factor for
all $a \in \A \setminus \{0,z\}$ and for all $0 \leq k < m$.

Moreover, put
\begin{equation}\label{eq:k_0}
    k_0 =   \begin{cases}
                -1 & \text{if $t \neq t_1 - 1$,} \\
                0 & \text{if $t = t_1 - 1$ and $t_2 \neq t_{m+1}$} \\
                \max\{\ell  \in \N \mid t_{\ell +1} \neq t_{m \oplus \ell} \} & \text{otherwise,}
            \end{cases}
\end{equation}
then the $k$-th max-$f$-image of $0^{t_1-1}$ is also
$(g_L^k(0,z))$-maximal factor for all $k_0 < k < m$.
\end{coro}

\pf

It holds
\begin{equation*}
    \Rex{00^{t_1 - 1}} = \{ k \in \A \setminus \{0\} \mid t_k = t_1 \text{ or } k = m \text{ and } t_{m+p} = t_1\}
\end{equation*}
and for all $a \in \A \setminus \{0\}$ we have $k \in \Rex{a0^{t_1
- 1}}$ if and only if $k = 0$ or the both following conditions are
satisfied:
\begin{itemize}
    \item[(i)] $z(k) = a - 1$ or $y(k) = a-1$,
    \item[(ii)] $t_k = t_1 - 1$ or $k = m$ and $t_{m+p} = t_1-1$.
\end{itemize}
The intersection of $\Rex{00^{t_1 - 1}}$ and $\Rex{a0^{t_1 - 1}}$
is not empty if and only if $a = z$ and $t = t_1 - 1$, in other
words, if and only if $0^{t_1 - 1}$ is a prefix of $0^tm$ what is
a LS factor having just two left extensions $0$ and $z$.

Similarly, we can prove that the $k$-th max-$f$-image of
$0^{t_1-1}$ is $(g_L^k(0,a))$-maximal factor for all $a \in \A
\setminus \{0,z\}$. Also similarly, the $k$-th max-$f$-image of
$0^{t_1-1}$, namely
$$
    \vpb^{k}(0^{t_1-1})\vpb^{k}(1)(k+1)^{-1},
$$
is $(g_L^k(0,z))$-maximal if it is not prefix of the LS factor
$$
\vpb^k(0^tm) = \vpb^k(0^t)\vpb^{k}(m)
$$
having the left extensions $g_L^k(0,z)$. The proof then follows
from Lemma~\ref{lem:simple_observations} part $(i)$ and
Lemma~\ref{lem:main_lemma} applied on $\vpb^{k}(1)$ and
$\vpb^{k}(m)$.

\pfk

Taking into account Lemmas~\ref{lem:lext_of_letters}
and~\ref{lem:generators_of_maximals}, one can prove the following
corollary using analogous techniques as in the proof of the
previous one. Note that $\Rex{\ell_0 0} = \{k \in \A \mid z(k-1) =
\ell_0-1 \text{ or } y(k-1) = \ell_0-1\}$ and $\Rex{a0} = \{1\}$
for all $a > \ell_0$, i.e., $0$ is $(\ell_0,\ell_0+z)$-maximal if
it is not a prefix of the $\ell_0$-th max-$f$-image of the
factor~\eqref{eq:the_max_factor} which reads
$$
\vpb^{\ell _0}(m)\vpb^{m + \ell_0}(1)(1 \oplus (m + \ell_0))^{-1}
= \vpb^{\ell _0}(m)\vpb^{m}(\ell_0+1)(1 \oplus (m + \ell_0))^{-1}.
$$

\begin{coro} \label{prop:maximals_t_1_eq_1}
If $t_1 = 1$, then the $k$-th max-$f$-image of $0$ is
$(g_L^k(\ell_0,a + \ell_0))$-maximal factor for all letters $a >
\ell_0$, $a \neq z$ and for all $0 \leq k < m - \ell_0$.

Moreover, $k$-th max-$f$-image of $0$ is $(g_L^k(\ell_0,z +
\ell_0))$-maximal if $k_0 \geq \ell_0$ and $k = k_0 - \ell_0, k_0
- \ell_0 +1,\ldots m - \ell_0$, where $k_0$ is defined
by~\eqref{eq:k_0}.
\end{coro}

\section{Affine complexity}

The aim of the present section is to find the necessary and
sufficient condition for the factor complexity of $\ubeta$ being
affine. In order the complexity to be affine, the first difference
of complexity $\Cd(n)$ must be constant. The following lemma says
when $\Cd(n)$ can change its value. The proof is an immediate
consequence of~\eqref{eq:connection}.
\begin{lem}
Let $\ub$ be an infinite word over a finite alphabet.
\begin{itemize}
    \item[(i)] If $\Cd(n+1) > \Cd(n)$, then the number of LS
    factor of length $n+1$ is greater then the number of LS
    factor of length $n$.
    \item[(ii)] If $\Cd(n+1) < \Cd(n)$, then $\ub$ contains
    $(a,b)$-maximal factor of length $n$ for some letters $a$ and
    $b$.
\end{itemize}
\end{lem}
That is, the complexity is affine if either $\ub$ does not contain
any $(a,b)$-maximal factor and all infinite LS branches have empty
common prefix or if each $(a,b)$-maximal factor of length $n$ is
``compensated" by appearing of a ``new'' LS factor of length
$n+1$. Examples of the first case are Arnoux-Rauzy words whose all
LS factors are prefixes of unique infinite LS branch. As for the
latter case, appearing of a ``new'' LS factor of length $n+1$
means there are a LS factor $v$ of length $n$ and its right
extensions $c$ and $d$ such that $vc$ and $vd$ are both LS, i.e
$v$ is the longest common prefix of two different LS factors --
Cassaigne~\cite{Cassaigne1997} call such LS factors \emph{strong
bispecial}.

Since $\ubeta$ comprises always at least one $(a,b)$-maximal
factor, each such $(a,b)$-maximal must be as long as the longest
common prefix of two different LS factors in order that the
complexity may be affine. We will prove that it is possible only
if the number of $(a,b)$-maximal factors is finite, thus in the
case of $d_\beta(1) = t_1(0\cdots0(t_1 - 1))^\omega$.
\begin{lem} \label{lem:not_affine_k_0}
    If $k_0 < m - 1$, where $k_0$ is defined by~\eqref{eq:k_0}, then the
    factor complexity of $\ubeta$ is not affine.
\end{lem}

\pf

If $k_0 < m - 1$, then the $(k_0 + 1)$-th max-$f$-image, if $t_1 >
1$, (resp. $(k_0 - \ell_0 + 1)$-th if $t_1 = 1$)  of $0^{t_1 -1}$
(resp. $0$) is $g_L^{k_0}(0,z)$-maximal. Consider the longest
common prefix of the LS factor $\vpb^{k_0}(0^tm)$ having left
extensions $g_L^{k_0}(0,z)$ and of the infinite LS branch
$\ubeta$, if $p>1$, or of the LS factor $\vpb^{m-1}(0^{t_1-1})$
with left extensions $m-1$ and $m$, if $p = 1$ (and so $\ubeta$ i
not an infinite LS branch). This factor equals
$\vpb^{k_0}(0^tm)(m\oplus k_0)^{-1}$ which is a prefix of the
$(k_0 + 1)$-th max-$f$-image (resp. $(k_0 - \ell_0 + 1)$-th if
$t_1 = 1$) of $0^{t_1 -1}$ (resp. $0$) and hence it is not
$(a,b)$-maximal for any distinct $a,b \in \A$. Overall, $\Cd(n_0)
< \Cd(n_0+1)$, where $n_0$ is the length of the factor
$\vpb^{k_0}(0^tm)(m\oplus k_0)^{-1}$.

\pfk

\begin{lem}
    If $d_\beta(1) = t_1(0\cdots0(t_1 - 1))^\omega$, then the factor complexity
    of $\ubeta$ is affine, namely $\C(n) = pn + 1, n \in \N$.
\end{lem}

\pf

In this case, $t = t_1 - 1$ and so $k_0 = 0 = m - 1$. Hence, the
$(0,a)$-maximal factor $0^{t_1-1}$ is at the same time the longest
common prefix of the only infinite LS branches $\ubeta$ and
$0^tm\vpb(0^tm)\vpb^2(0^tm)\cdots$. But $0^{t_1-1}$ is the only
$(a,b)$-maximal and prefixes of these two infinite LS branches are
the only LS factors of $\ubeta$, thus, the proof is complete.

\pfk

\begin{lem}
If $\beta \in \S$ and $d_\beta(1) \neq t_1(0\cdots0(t_1 -
1))^\omega$, then the factor complexity of $\ubeta$ is not affine.
\end{lem}
\pf

In the case when $p>1$, there are $m+1$ infinite LS branches given
by Proposition~\ref{prop:list_of_inf_branches}. Let us denote them
by $\ub_0 = \ubeta, \ub_1, \ldots, \ub_m$ and put
$$
    n_0 = \max\{|v| \mid v = \lcp{\ub_i}{\ub_j}, i \neq j, i,j =
    0,1,\ldots,m\}.
$$
We have $\Cd(n) \geq \#\Le{\ub_0} - 1 + \sum_{k=1}^m \#\Le{\ub_k}
-1 \geq p-1+m$ for all $n > n_0$. Due to
Corollary~\ref{prop:maximals_0t_m}, we know that there exist
infinitely many $(g_L^k(0,z))$-maximal factors, $k = 0,1,\ldots$,
and hence there must exist a LS factor of length $n_1 > n_0$ which
is not a prefix of any LS branch and so $\Cd(n_1) > m + p - 1 =
\Cd(1)$.

In the case of $p = 1$, the proof is analogous. Only difference is
that there are only $m$ infinite LS branches since $\ubeta$ is
not. \pfk

\begin{rmrk} For the word $\ubeta$ with $d_\beta(1) = t_1(0\cdots0(t_1 - 1))^\omega$
 we may easily describe all left special factors.
 If the length of the period $p > 1$,
each LS factor is a prefix of one of two infinite LS branches
$\ubeta$ and $0^{-1}\ubeta$. If $p=1$, then $\ubeta$ is not an
infinite LS branch and so each LS factor is prefix of the unique
infinite LS branch $0^{-1}\ubeta$. Hence, we obtain the known
result that $\ubeta$ is Sturmian if and only if $d_\beta(1) =
t_1(t_1 - 1)^\omega$. We were pointed out by Christiane Frougny
that numbers $\beta$ satisfying $d_\beta(1) = t_1(0\cdots0(t_1 -
1))^\omega$ are Pisot units. Such Parry number $\beta$ is a root
of the polynomial $x^{p+1}- t_1x^p -x+1$.
\end{rmrk}

\begin{lem}
    Let $\beta \notin \S$ and let $k_0 \geq m-1$. Then the
    factor complexity of $\ubeta$ is not affine.
\end{lem}

\pf

As shown in the proof of Lemma~\ref{lem:not_affine_k_0}, the
$k$-th max-$f$-image of $0^{t_1-1}$ (resp. $0$ if $t_1 = 1$), $k =
0,1,\ldots,m-1$, is not beginning in $0^tm$ and it is equal to the
longest prefixes of some two LS factors. In order that the
complexity is affine, also all max-$f$-images of the
factor~\eqref{eq:the_max_factor} must be as long as the longest
prefixes of some two LS factors.

Let $t_1 > 1$. Then the factor~\eqref{eq:the_max_factor} must be
of the same length as the longest common prefix of $\ubeta$ and
$m$-th max-$f$-image of its own -- remember that the longest
common prefix of $\ubeta$ and $k$-th max-$f$-image
of~\eqref{eq:the_max_factor} is the $k$-th max-$f$-image
of~$0^{t_1-1}$ for $k = 0,1,\ldots,m-1$. Formally,

\begin{multline*}
    |0^tm\vpb^m(0^{t_1-1}1)(1+m)^{-1}| = \\
    = |\lcp{\ubeta}{\vpb^m(0^tm)\vpb^{2m}(0^{t_1-1}1)(1\oplus (2m))^{-1}}|=
    |\vpb^m(0^t m)(m\oplus m)^{-1}|
\end{multline*}
which is never satisfied for $|\vpb^m(0^t m)| \leq
|\vpb^m(0^{t_1-1} 1)|$.

Let $t_1 = 1$. Following the same reasoning as for the case $t_1
> 1$, a necessary condition for the complexity to be affine is
that the factor~\eqref{eq:the_max_factor}
$$
    m\vpb^m(1)(1+m)^{-1}
$$
must be of the same length as the longest common prefix of the
$(m-\ell_0)$-th max-$f$-image of its own and $\ubeta$, namely
$$
    |\lcp{\ubeta}{\vpb^{m-\ell_0}(m)\vpb^{2m-\ell_0}(1)(1\oplus
    (2m-\ell_0))^{-1}}|
    = |\vpb^{m-\ell_0}(m)(m\oplus (m-\ell_0))^{-1}|
$$
which is never satisfied for $|\vpb^{m-\ell_0}(m)| \leq
|\vpb^m(1)|$. \pfk

Putting all lemmas of this section together, we obtain the main
theorem of this paper.
\begin{thm} \label{thm:the_main_one}
    Let $\beta$ be a non-simple Parry number. The factor
    complexity of $\ubeta$ is affine if and only if $d_\beta(1) = t_1(0\cdots0(t_1 -
    1))^\omega$.
\end{thm}

\section{Conclusion}
Among infinite words $\ubeta$ associated with Parry numbers we may
identify Arnoux-Rauzy words. An infinite word is said to be
Arnoux-Rauzy of order $\ell$, if for any length $n\in \mathbb{N}$
there exists exactly one left special factor and one right special
factor both of length $n$ and, moreover, these special factors
have just $\ell$ left and $\ell$ right extensions respectively.
Arnoux-Rauzy words can be considered as a natural generalization
of Sturmian words to more letter alphabets.

Is is easy to see that only Sturmian words among $\ubeta$
correspond to $\beta$ with $d_\beta(1) = t_11$ or $d_\beta(1) =
t_1(t_1-1)^\omega$. The word $\ubeta$  is an Arnoux-Rauzy word of
order $m\geq 3$ if and only if $d_\beta(1) = t_1^{m-1}1$. It means
that there is no Arnoux-Rauzy word over more letter alphabet
associated with non-simple Parry number. A direct consequence of
the definition of Arnoux-Rauzy words is that the complexity of
Arnoux-Rauzy word is affine and that any left (right) special
factor is a prefix of an infinite left (right) special branch.

In the previous section, we have proved that the infinite word
$\ubeta$ associated with a non-simple Parry number $\beta$ has the
affine complexity if and only if $d_\beta(1) = t_1(0\cdots0(t_1 -
1))^\omega$. In fact, we have proved that the complexity is affine
if and only if any left special factor of $\ubeta$ is a prefix of
an infinite left special branch.  The validity of the same
statement for infinite words associated with simple Parry numbers
is proven in~\cite{Bernat2007}. However, this equivalency is not a
general rule for the factor complexity of fixed points of
primitive morphisms. For a counter example see~\cite{Chacon1969}
and~\cite{Ferenczi1995}.

It is known that Sturmian words have many equivalent definitions,
see~\cite{Berstel2007a} for more. In 2001
Vuillon~\cite{Vuillon2001} showed that a binary infinite word is
Sturmian if and only if each its factor has exactly two return
words. In the article~\cite{Vuillon2000} Vuillon introduced the
property $R_\ell$: an infinite word satisfies the property
$R_\ell$ if each its factor has exactly $\ell$ return words.
Therefore, words with $R_\ell$ can be considered as another
generalization of Sturmian words. In~\cite{Justin2000} Justin and
Vuillon proved that Arnoux-Rauzy words of order $\ell$ have the
property $R_\ell$. Applying Theorem~4.5 of~\cite{Balkova2008}, we
see that all $u_\beta$ with affine complexity have also the
property $R_\ell$.

\section*{Acknowledgement}

We thank Christiane Frougny for many fruitful discussions. We
acknowledge financial support by the grants MSM6840770039 and
LC06002  of the Ministry of Education, Youth, and Sports of the
Czech Republic.

\bibliographystyle{plain}
\bibliography{kloudak_bibtex}

\end{document}